\input amstex.tex
\input amsppt.sty   
\magnification 1200
\vsize = 8.28 true in
\hsize=6.2 true in
\nologo
\NoRunningHeads        
\parskip=\medskipamount
        \lineskip=2pt\baselineskip=18pt\lineskiplimit=0pt
       
        \TagsOnRight
        \NoBlackBoxes

        \topmatter
        \title
        Space Quasi-Periodic Standing Waves for
        \\Nonlinear Schr\"odinger Equations
        \endtitle
\author
         W.-M.~Wang \endauthor     
\address
{CNRS and D\'epartement de Math\'ematique, Universit\'e Cergy-Pontoise, 95302 Cergy-Pontoise Cedex, France}
\endaddress
        \email
{wei-min.wang\@math.cnrs.fr}
\endemail
\abstract
We construct {\it space} quasi-periodic standing wave solutions to 
the nonlinear Schr\"odinger equations on $\Bbb R^d$ for
{\it arbitrary} $d$. This is a type of quasi-periodic nonlinear Bloch-Floquet waves. 
\endabstract
    
        \bigskip\bigskip
        \bigskip
        \toc
        \bigskip
        \bigskip 
        \widestnumber\head {Table of Contents}
        \head 1. Introduction to the Theorem
        \endhead
        \head 2. Green's function estimates in $(\theta_1, \theta_2, ..., \theta_d)$
                \endhead
        \head 3. Nonlinear analysis -- proof of the Theorem\endhead 
        \endtoc
        \endtopmatter
        \vfill\eject
        \bigskip
       
\document
\head{\bf 1. Introduction to the Theorem}\endhead
Consider the nonlinear Schr\"odinger equations (NLS) on $\Bbb R^d$:
$$
i\frac\partial{\partial t}U =-\Delta U-|U|^{2p}U,\tag 1.1
$$
where $p\geq 1$ and $p\in\Bbb N$ is arbitrary; $U$ is a complex valued function on $\Bbb R\times \Bbb R^d$.
In this paper, we seek standing wave solutions of the form $$U(t, x)=e^{-iEt}u(x),\tag 1.2$$
where $E\in\Bbb R$, and $u$ is {\it even} and {\it quasi-periodic} in each $x_k$, $k=1, 2, ..., d$, given by a quasi-periodic cosine
series: 
$$u(x)=u(x_1, x_2, ..., x_d)=\sum_{j_1, j_2, ..., j_d} \hat u(j_1, j_2, ..., j_d)\prod_{k=1}^d\cos(j_k\cdot\lambda_k)x_k,\tag QP$$
where for each $k\in\{1, 2, ..., d\}$, $j_k\in\Bbb Z^2$ and $\lambda_k\in (1/2, 3/2)^2$. The $\lambda_k$'s
are the {\it parameters} in the problem, and are assumed to be irrational, satisfying 
$$\Vert j_k\cdot\lambda_k\Vert_{\Bbb T}\neq 0,\tag D$$ for all $j_k\neq 0$,
where $\Vert\cdot\Vert_{\Bbb T}$ denotes distance to the integers.
We note that the quasi-periodic series (QP) reduces to a periodic cosine series if $j_k$ and $\lambda_k$ were
one dimensional: $j_k\in\Bbb Z$ and
$\lambda_k\in(1/2, 3/2)$ for $k=1, 2, ..., d$.  For example, setting $\lambda_k=1$ for all $k$, leads to a periodic
series with period $2\pi$ in each directions. In that case, solutions with more general time dependence, the time quasi-periodic solutions,
are known to exist from e.g., \cite{W1}, cf. also \cite{W2}.  

Substituting the Ansatz (1.2) into (1.1) yields the following stationary, nonlinear elliptic problem: 
$$-\Delta u-|u|^{2p}u=Eu.\tag 1.3$$
For $u\in H^1(\Bbb R^d)$ with a {\it fixed} $L^2(\Bbb R^d)$ norm, there is a well established variational structure under appropriate 
conditions on $p$: $E$
is a Lagrange multiplier and (1.3) are the minimizers
for the energy functional:
$$\Cal E(U)=  \int_{\Bbb R^d} dx[\frac{1}{2}\Vert \nabla U\Vert^2 -\frac{1}{2p+2} |U|^{2p+2}].$$
Localized standing wave solutions are well known from the works of, for example, Cazenave and Lions \cite{CL}. 
(Cf. also the references therein.) The $u$'s given by (QP), even though smooth,
are only in $L^\infty$ and solving (1.3) produces space quasi-periodic nonlinear Bloch-Floquet waves, which are   
{\it not} localized. (For quasi-periodic {\it linear} Bloch-Floquet theory in one dimension, see e.g.,  \cite{DS, E}, cf. also \cite{K} for {\it linear} quasi-periodic ground states in arbitrary dimensions.)

\noindent{\it Remark.} For the purpose of this paper, the sign of the nonlinear term is 
unimportant, i.e., it can be focusing or defocusing. (See the remark after the Theorem.)
Functions that are even under $x_k\to -x_k$, for all $k=1, 2, ..., d$, form an invariant subspace 
for (1.3). Here we seek solutions $u$ in this subspace given by the series in (QP). 

To simplify notations, define 
$$ (j\cdot\lambda)^2:=\sum_{k=1}^d(j_k\cdot\lambda_k)^2,\tag 1.4$$
where $j=(j_1, j_2, ..., j_d)\in\Bbb Z^{2d}$ and $\lambda=(\lambda_1, \lambda_2, ..., \lambda_d)\in (1/2, 3/2)^{2d}$.  

Let $$\tilde U=ae^{-i{({\tilde j}\cdot\lambda)}^2 t}\prod_{k=1}^d\cos(\tilde j_k\cdot\lambda_k)x_k,\tag 1.5$$
where $a\in\Bbb R$ and ${\tilde j}\in \Bbb Z^{2d}$. (If $\tilde j=0$, $\tilde U= ae^{i|a|^{2p}t}$ trivially solves (1.1).) 
Then $\tilde U$ satisfies the linear equation:
$$
i\frac\partial{\partial t}\tilde U=-\Delta \tilde U;\tag 1.6
$$
and $$\tilde u=a\prod_{k=1}^d\cos(\tilde j_k\cdot\lambda_k)x_k\tag 1.7$$
satisfies
$$
-\Delta \tilde u=\tilde E \tilde u,\tag 1.8
$$
with $$\tilde E=(\tilde j\cdot\lambda)^2.\tag 1.9$$

Our main result is 

\proclaim{Theorem} For every solution to the linear equation (1.6) in the form (1.5)
$$\tilde U=ae^{-i{({\tilde j}\cdot\lambda)}^2 t}\prod_{k=1}^d\cos(\tilde j_k\cdot\lambda_k)x_k,$$
where $a\in\Bbb R$ and ${\tilde j}\in \Bbb Z^{2d}$, there is a set in $\lambda$, $\Lambda\subset (1/2, 3/2)^{2d}$
satisfying $$\text{meas } \Lambda\geq 1-|a|^{p/6},\tag 1.10$$ 
provided $|a|\ll 1$ . If $\lambda\in\Lambda$, then there is
a solution $U$, bifurcating from $\tilde U$, to the nonlinear equation (1.1) in the form (1.2, QP):  
$$U(t, x)=e^{-i{[({\tilde j}\cdot\lambda)}^2+\Cal O(|a|^{2p})] t} [a\prod_{k=1}^d\cos(\tilde j_k\cdot\lambda_k)x_k+ \Cal O(|a|^{p})].\tag 1.11$$
The nonlinear eigenvalue $E$, as a function in $\lambda$, is $C^1$ on $(1/2, 3/2)^{2d}$.
\endproclaim
\noindent{\it Remark.} For notational simplicity we have taken $\lambda_i$, $i=1, 2, ..., d$, to be two-dimensional. The 
Theorem holds for higher dimensional $\lambda_i$, $i=1, 2, ..., d$, with essentially the same proof. The set $\Lambda$ is
a Cantor set (of positive measure). Since $a$ is small, the same Theorem holds if the nonlinearity enters with a plus sign (defocusing). 

\subheading{1.1. Some background}

Most of the results in the literature on (1.1) or (1.3) are for $u$, which are fast decaying or periodic in $\Bbb R^d$. 
To our knowledge, the above Theorem is the first such result on global in time, non-decaying solutions $u$ 
which do not have an underlying translation symmetry group. 
(Cf. Moser \cite{M} for an iterative method in the space periodic setting, i.e., on the quotient space
$L^2(\Bbb R^d/\Bbb Z^d):= L^2(\Bbb T^d)$.)
It is periodic in time (with only the basic frequency), quasi-periodic in space and exists in arbitrary dimensions. 
The Theorem shows that under appropriate conditions, every small even generalized eigenfunctions of the linear operator in (1.8) bifurcates
to an eigenfunction of the nonlinear operator in (1.3), after small deformation. 

Generally speaking, due to the non-compact $\Bbb R^d$ setting, there are very few known results on space quasi-periodic solutions to 
nonlinear partial differential equations. In one dimension, Damanik and Goldstein proved the global existence and uniqueness to Cauchy problems for 
the KdV equation with small quasi-periodic initial data \cite{DG}. Their method, however, seems to hinge on the integrable structure.
It is noteworthy that the Cauchy solutions in \cite{DG} are almost-periodic in time (and quasi-periodic in space). This result in fact 
motivated us to seek space quasi-periodic solutions in a more general setting, albeit with simpler time dependence,
as in the Theorem. However, solutions with more complicated time dependence to the NLS in (1.1), such as space-time quasi-periodic solutions 
can be analyzed, see \cite {W3}. 
Note also that equation (1.1) is used to study Bose-Einstein
condensation, cf. e.g., \cite{LOSK}, and is usually called the Gross-Pitaevskii equation, when 
seeking non-decaying solutions. 

One may pose similar questions for nonlinear difference equations, for example, for the
Frenkel-Kontorova model, studied in Aubry-Mather theory, cf. e.g., \cite{EFRJ} for its physical origin and \cite {SdlL,~GPT}
for KAM-type results in one dimension. The method proposed here should be
applicable, providing (time periodic) space quasi-periodic solutions in arbitrary dimensions, corresponding to {\it sliding}.

\subheading{1.2. Ideas of the proof} 

Since $\tilde u$ is real, we may seek real solutions $u$ to (1.3). 
Use diag $\cdot$  to denote a diagonal matrix. Substituting 
(QP) into (1.3) leads to the nonlinear matrix equation on $\ell^2(\Bbb Z^{2d})$: 
$$\text{diag }(\sum_{k=1}^d(j_k\cdot\lambda_k)^2-E)\hat u-(\hat u)^{*2p}*\hat u=0;$$
with the linearized operator being
$$H=\text{diag }(\sum_{k=1}^d(j_k\cdot\lambda_k)^2-E)-(2p+1)(\hat u)^{*2p}*.$$
To fix ideas, set $\hat u ={\hat u}^{(0)}$. $H$ is then quasi-periodic in $d$-dimensions
on diagonal plus a convolution operator. The issue is to control the inverse of $H$.
The main difficulty here is that for $d>1$, Diophantine conditions on $\lambda$ do not suffice. The problem is more geometric,
and we use the semi-algebraic set technique developed by Bourgain in the study of Anderson localization \cite{B3} to
do the linear analysis. (Such techniques first appeared in \cite{BGS} on quasi-periodic Anderson localization in $\Bbb Z^2$.) 
This is different from the space periodic setting in \cite{W1}, cf. also \cite{W2}, where 
the quasi-periodicity is in time only, which is one dimensional. Diophantine conditions together with eigenvalue variations suffice for the
linear analysis. (The main work in \cite{W1,~2} is to {\it extract parameters} from the nonlinear term, 
in order to deal with the {\it original (fixed)} nonlinear equations such as that in (1.1), and {\it not} merely a 
family of parameter-dependent equations as in \cite{B1, 2}. The method of extraction is algebraic.)

Once we have good control on the inverse of $H$, the nonlinear analysis proceeds
using a Newton iteration, based on Chap.~18 of \cite{B2}, cf. \cite{BW, W1, 2}. This
part is rather standard, and shares many common features with other KAM-type schemes.
It is, in fact, simpler here, since the ``dynamical variables" are 
the space variables $j\in\Bbb Z^{2d}$ and there is {\it no} modulation to the frequency $\lambda$. 
\bigskip
\noindent{\bf Acknowledgement.} It is a pleasure to thank T. Spencer and R. de la Llave for discussions.

\head{\bf 2. Green's function estimates in $(\theta_1, \theta_2, ..., \theta_d)$}\endhead
Returning to the problem at hand, we seek solutions $U$ close to $\tilde U$, in the form (1.2) and (QP),
leading to the nonlinear matrix equation on $\ell^2(\Bbb Z^{2d})$: 
$$\text{diag }(\sum_{k=1}^d(j_k\cdot\lambda_k)^2-E)\hat u-(\hat u)^{*2p}*\hat u=0, \tag 2.1$$
and the linearized operator 
$$H=\text{diag }(\sum_{k=1}^d(j_k\cdot\lambda_k)^2-E)-(2p+1)(\hat u)^{*2p}*.$$

To fix ideas, we give an example of the convolution matrix. 
Let $d=2$ and $p=1$; write $(x, y)$ for $(x_1, y_1)$ and set $\hat u=\hat{\tilde u}$, as in (1.7); set 
$a=1$. We compute:
$$\aligned &\cos^2 (\tilde j_1\cdot\lambda_1)x \cos(j_1\cdot\lambda_1)x\\
=&1/2[\cos( j_1\cdot\lambda_1)x]+1/4[\cos(j_1+2\tilde j_1)\cdot\lambda_1x+\cos(j_1-2\tilde j_1)\cdot\lambda_1x];\endaligned$$
and similarly in the $y$ variable, i.e., with the subindex $1\leftrightarrow 2$ and $x\leftrightarrow y$.

It follows that $\hat{\tilde u}^{*2}$ is a convolution matrix on $\Bbb Z^4$, with the non-vanishing entries:
\item {$\bullet$} $1/4$ on diagonal; 
\item {$\bullet$} $1/8$ on the $(\pm 2\tilde j_1, 0)$ and $(0, \pm 2\tilde j_2)$ off-diagonals ($4$ in total); 
and
\item {$\bullet$} $1/16$ on the $(\pm 2\tilde j_1, \pm 2\tilde j_2)$ off-diagonals ($4$ in total).
\bigskip
We now proceed to the analysis.
\smallskip
\subheading {2.1. Lyapunov-Schmidt decomposition} 

We use a Newton scheme to solve (2.1), using as initial approximation $u^{(0)}=\tilde u$ in (1.7) and $E^{(0)}=\tilde E=(\tilde j\cdot\lambda)^2$ in (1.9).
In the matrix notation of (2.1), $\hat u$ is a column vector and ${\hat u}^{(0)}$ is the column vector with 
${\hat u}^{(0)}(j)=a/2^d$, if $j_k=\pm\tilde j_k$, $k=1, 2, ..., d$, and $0$ otherwise. Below, since we only work with $\hat u$, we slightly abuse the notation and 
write $u$ for $\hat u$. We may also assume $a>0$ as, if $u$ is a solution, then so is $-u$.

Let $$\Cal S=\{\pm \tilde j_k, k=1, 2, ..., d\}.\tag 2.2$$ Writing (2.1) as $$F(u)=0,\tag 2.3$$ the equations are divided into the 
$Q$-equations: 
$$F(u)|_{\Cal S}=0;\tag 2.4$$
and the $P$-equations:
$$F(u)|_{\Bbb Z^{2d}\backslash\Cal S}=0.\tag 2.5$$
The amplitudes on the set $\Cal S$ are held fixed: $$u|_\Cal S=a/2^d;\tag 2.6$$ while the $Q$-equations are used to solve for $E$. Due to symmetry, the $2^d$ equations in (2.4) are the same, yielding
$$E=(\tilde j\cdot\lambda)^2-(2^d/a)(u)^{*2p+1}|_{\tilde j}.\tag 2.7$$
So, for example, the first iteration gives 
$$E^{(1)}=(\tilde j\cdot\lambda)^2-(2^d/a)(u^{(0)})^{*2p+1}|_{\tilde j}.\tag 2.8$$

Substituting the result in (2.7) into the $P$-equations (2.5), we use a Newton scheme to solve for $u$
on ${\Bbb Z^{2d}\backslash\Cal S}$. For simplicity, omitting the subindex ${\Bbb Z^{2d}\backslash\Cal S}$ from now on, 
we have formally, (note that the $\Delta$ below denotes increment), 
$$\Delta u= -[F'(u)]^{-1}F(u),\tag 2.9$$
where $F'(u)$ is the linearized operator: 
$$F'(u)=\text{diag }(\sum_{k=1}^d(j_k\cdot\lambda_k)^2-E)-(2p+1)(u)^{*2p}*.\tag 2.10$$

Generally speaking, the idea is to start with the initial approximation 
$(u^{(0)}, E^{(0)})$ as in (1.7, 1.9) and to 
iterate the Newton scheme, with each iteration $i$ resulting in 
an approximate solution $(u^{(i)}, E^{(i)})$, after appropriate excisions in $\lambda$; and as $i\to\infty$, $(u^{(i)}, E^{(i)})$ 
converges to a solution $(u, E)$ to (1.3). Hence $U$ in (1.2)-(QP) is a solution to (1.1) for a subset of $\lambda$. 

\noindent {\it Remark.} We note that the above $P$ and $Q$-equations are decomposed
according to the Fourier support of $\tilde u$, $\Cal S$, and uses the condition (D).
The $Q$-equations are resonant, as the diag in (2.10)  is $0$ on $\Cal S$,
when $E=\tilde E$; while the $P$-equations are non-resonant.  
 
\subheading {2.2. Invertibility of the linearized operators}  

From (2.9), the invertibility of $F'$ is central to the Newton iteration.  
Since we seek solutions close to $u^{(0)}$, which is only supported on $\Cal S$, 
we adopt a {\it multiscale} Newton scheme. The idea is as follows. 

At each iteration $i$, choose an appropriate scale $N$ and estimate $[F'_N]^{-1}$, 
where $F'_N$ is $F'$ restricted to $$[-N, N]^{2d}\subset \Bbb Z^{2d}.\tag 2.11$$ 
We call the $[F'_N]^{-1}$, the Green's functions. 
To facilitate the estimates, add $d$ auxiliary variables 
$$\theta_1, \theta_2, ..., \theta_d,$$
to $F'$ and define:
$$F'(\theta_1, \theta_2, ..., \theta_d):=\text{diag }(\sum_{k=1}^d(j_k\cdot\lambda_k+\theta_k)^2-E)-(2p+1)(u)^{*2p}*.\tag 2.12$$
Denote $(\theta_1, \theta_2, ..., \theta_d)$ by $\theta\in\Bbb R^d$.
We first make estimates on $F'_N(\theta)$ in $\theta$ and then use the covariance with respect to the $\Bbb Z^{2d}$ action on $\Bbb R^d$: 
$$(\theta_1, \theta_2, ..., \theta_d)\mapsto (\theta_1+j_1\cdot \lambda_1, \theta_2+j_2\cdot \lambda_2, ..., \theta_d+j_d\cdot \lambda_d),\tag 2.13$$ 
to deduce estimates for $$[F'_N(\theta=0)]^{-1}:=[F'_N]^{-1},$$ the Green's functions used in the Newton scheme (2.9).

\subheading {2.3. The $(\theta_1, \theta_2, ..., \theta_d)$ estimates}

Denote the linearized operator $F'$ by $T$; and $F'_N$, $T_N$. The goal of this
section is to estimate the Green's functions $T_N^{-1}(\theta)$ for all $N$, away from a set in $\theta\in\Bbb R^d$ of small sectional measure, after {\it appropriate} excisions in $\lambda\in (1/2, 3/2)^{2d}$. 
To apply the covariance in (2.13) for the Green's function analysis, it is essential that the excised set in $\lambda$ is {\it independent} of the starting point, i.e., at $j_1=j_2=...=j_d=0$, of $\theta\in\Bbb R^d$, 
which we denote by $\vartheta\in\Bbb R^d$.
This is, in essence, accomplished by variable reduction, eliminating the variable $\vartheta\in\Bbb R^d$. (See Lemmas~2.2 and 2.3 below). 

Since $d$ is arbitrary, the geometry of the sets in $\theta=(\theta_1, \theta_2, ..., \theta_d)$ comes into play. Diophantine conditions, i.e., quantitative versions of (D), generally do not suffice, and we shall 
use the semi-algebraic set technique developed by Bourgain in \cite {B3}, cf. Chap.~9 \cite{B2}. For that purpose, we need that $u^{(i)}$ and $E^{(i)}$ are algebraic in $\lambda$
and control their degrees. To begin with, $u^{(0)}$ does not depend on $\lambda$ (recall that $u^{(0)}$ now stands for ${\hat u}^{(0)}$), and from (1.9), (2.8), $E^{(0)}$ and $E^{(1)}$
are both quadratic polynomials in $\lambda$. 

Since $u^{(i)}$ and $E^{(i)}$ depend on the 
scale $N$, we denote them by $u_N$ and $E_N$ in this section.
We assume what is needed on $u_N$ and $E_N$ in (2.10) from the nonlinear analysis, in order to estimate the Green's functions. 
Later in sect.~3, we verify these assumptions.  

Let us first define a semi-algebraic set. 

\noindent{\bf Definition.} A set $S$ is called semi-algebraic if it is a finite union of sets defined by a 
finite number of polynomial equalities and inequalities. More specifically, let $\Cal P=\{P_1, P_2, ..., P_s\}\subset \Bbb R[x_1, x_2, ..., x_n]$
be a family of $s$ real polynomials of degree bounded by $\kappa$. A (closed) semi-algebraic set $S$ is given by an expression
$$S=\bigcup_j \bigcap_{\ell\in\Cal L_j} \{P_\ell s_{jl}0\},\tag S$$
where $\Cal L_j\subset \{1, 2, ..., s\}$ and $s_{jl}\in\{\geq, =,\leq\}$ are arbitrary. We say that $S$ as introduced above has degree at most $s\kappa$ and
its degree $B$ is the minimum $s\kappa$ over all representations (S) of $S$. 

The following is a special case of Theorem 1 in \cite{Ba}, cf. Theorem 9.3 in Chap.~9 \cite{B2}.

\proclaim{Lemma~2.1}Let $S\subset \Bbb R^n$ be as in {\rm (S)}. Then the number of connected components of $S$ does not exceed 
$\Cal O(s\kappa)^n$. 
\endproclaim

The two properties of semi-algebraic sets that play a central role here are the Tarski-Seidenberg principle, which states that the projection 
of a semi-algebraic set of $\Bbb R^n$ onto $\Bbb R^{n-1}$ is semi-algebraic; and the Yomdin-Gromov triangulation theorem of these sets. They are both 
stated in \cite{B3}, cf. the references therein. (For the complete proof of the Yomdin-Gromov triangulation theorem, see \cite{BiN}, cf., also the earlier 
paper \cite{Bu}.) We do not repeat them here, except their consequences for thin sets. 

Below we call connected open sets intervals. Our main goal is to prove the following.

\proclaim{Main Lemma} Let $I$ be an interval in $(1/2, 3/2)^{2d}$, $u_N$ and $E_N$ two sequences of real rational functions in $\lambda$,
$$u_N:\, I\mapsto \ell^2(\Bbb Z^{2d}),$$
satisfying $$u_N(j)=0, \, j\notin[-N^K, N^K]^{2d},$$
for some $K>1$; 
and $$E_N: I\mapsto \Bbb R.$$
For $0<a\ll 1$, assume that there exists $N_0=N_0(a)\gg1$, such that for $N\geq N_0$,  the following conditions are satisfied:
$$\text{deg }u_N\lesssim e^{(\log N)^3},\tag 2.14$$
$$|u_N(j)|\leq e^{-\gamma|j|},\, j\in \Bbb Z^{2d} \, (\gamma>0),\tag 2.15$$
$$\Vert u_N-u_{N+1}\Vert_{\ell^2(\Bbb Z^{2d})} \leq  e^{-\tilde\gamma N}\,(\tilde\gamma>\gamma>0);\tag 2.16$$
and 
$$\text{deg }E_N\lesssim e^{(\log N)^3}.\tag 2.17$$
$$E_N=\Cal O(1)\tag 2.18$$
$$|E_N-E_{N+1}|\leq e^{-\tilde\gamma N},\tag 2.19$$
For all $N\geq N_0$, there exists $\Cal A_N\subset I$, a semi-algebraic set of 
$$\align \text{deg }&\Cal A_N\leq N^{8d},\tag 2.20\\
\text{meas }(\Cal A_{N-1}\backslash&\Cal A_N)\leq N^{-c},\, c>0,\tag 2.21\endalign$$
such that for any $\lambda\in \Cal A_N$, there exists a subset $\Theta_N\subset \Bbb R^d$, whose sectional measures
satisfy
$$\text{meas } [\theta_i| \forall \text{ fixed }\theta_k, k\neq i; \theta\in \Theta_N]\leq e^{-N^\tau} \, (\tau>0),\tag 2.22$$
for all $i=1, 2, ..., d$. If $\theta\notin\Theta_N$, the linearized operator $F':=T$ in (2.12), after truncations, satisfy the estimates
$$\Vert [T_N(u_N, E_N)(\theta)]^{-1}\Vert_{\text{Op}}\leq e^{N^\sigma}\, (1>\sigma>\tau>0),\tag 2.23$$
and
$$|[T_N({u_N, E_N})(\theta)]^{-1}(j,j')|\leq e^{-\beta |j-j'|} (0<\beta<\gamma),\, \forall |j-j'|>N/10.\tag 2.24$$
\endproclaim

The nonlinear construction in sect.~3 will verify (2.14)-(2.19) by using the
double exponential convergence of the Newton scheme. Note that the algebraic in $\lambda$
requirements on $u_N$ and $E_N$ for the {\it linear} analysis is due to quasi-periodicity in space, 
this is different from Chaps.~19 and 20 in \cite{B2} and \cite{W1,~2}, which
are space periodic and the algebraic dependence is only used in the nonlinear analysis.  

\subheading{2.4. Proof of the Main Lemma} 

The proof is an application of Proposition~2.2 in \cite{B3}, complemented as Theorem~4.1 in \cite{JLS}. 
The algebraic arguments rely on Lemmas~1.18 and 1.20 in \cite{B3}, and are stated below as Lemmas~2.2 and 2.3.
The analysis arguments after formula (2.42), p 696-699  in sect.~2 \cite {B3} have been dissected and clarified 
in Theorems~3.6, 4.1 and their proofs in \cite{JLS}. 

\proclaim{Lemma~2.2}
Let $A\subset [0, 1]^{n+r}$ be semi-algebraic of degree $B$ and such that  
$$ \text{for each }t\in[0, 1]^r, \, \text{meas}_nA(\cdot\,,t)<\eta, \eta>0.\tag 2.25$$
Then $$\Cal A:=\{(x_1, x_2, ..., x_{2^r})|A(x_1)\cap...\cap A(x_{2^r})\neq \emptyset\}\subset[0, 1]^{n2^r}\tag 2.26$$
is semi-algebraic of degree at most $B^C$ and measure at most 
$$\eta_r=B^C\eta^{n^{-r}2^{-\frac{r(r-1)}{2}}}\tag 2.27$$
with $C=C(r)>1$. 
\endproclaim

Lemma~2.2 is a variable reduction lemma, eliminating the $r$-dimensional variable $t$.
It is worth noting that $2^r$ copies of $A$ are used. The measure in (2.27), however, is in
$n2^r$ dimensions; while we need the measure of a $n$-dimensional section of $\Cal A$.
Lemma~1.20 in \cite{B3} serves this purpose, and is stated below.

\proclaim{Lemma 2.3} Let $A\subset [0, 1]^{n\rho}$ be a semi-algebraic set of degree $B$ and 
$$\text{meas}_{n\rho}A<\eta.$$
Let $\omega_i\in [0,1]$, $i=1, 2, ..., n$,
and $$\omega=(\omega_1, \omega_2, ..., \omega_n)\in [0, 1]^{n}.$$
Let $k_i\in \Bbb Z$, $i=1, 2, ..., n$,
and $$k=(k_1, k_2, ..., k_n)\in \Bbb Z^{n}.$$
Denote by $\{\cdot\}$, the fractional part of a real number in $[0, 1)$, and   
$$k\omega:=(\{k_1\omega_1\}, \{k_2\omega_2\}, ..., \{k_n\omega_n\}).\tag 2.28$$
Let $K_1$, $K_2$, ..., $K_{\rho-1}\subset \Bbb Z^{n}$ be finite sets with the following properties:
$$\min_{1\leq \ell \leq n}|k_\ell|>[B \max_{1\leq \ell' \leq n}|m_{\ell'}|]^C,\tag 2.29$$
if $k\in K_i$ and $m\in K_{i-1}$, $i=2, ..., \rho-1$, and where $C=C(n, \rho)$. Assume also 
$$\frac{1}{\eta}>\max_{k\in K_{\rho-1}}|k|^C.\tag 2.30$$
Then 
$$\aligned &\text{meas }\{\omega\in [0, 1]^{n}|(\omega, k^{(1)}\omega, ..., k^{(\rho-1)}\omega)\in A \text{ for some } k^{(i)}\in K_i\}\\
<&B^C\delta,\endaligned\tag 2.31$$
where $$\frac{1}{\delta}=\min_{k\in K_1}\min_{1\leq \ell\leq n}|k_\ell|.\tag 2.32$$
\endproclaim

{\noindent}{\it Remark.} As noted  in the first paragraph in sect.~2.3, Lemmas~2.2 and 2.3 are the tools to eliminate the variable $\vartheta$ (the starting point).
The inequalities in (2.29) are {\it steepness} conditions.
Even though we will not make use of it in proving the Main Lemma, we mention that 
using the special structure of the $\Bbb Z^{2d}$ action
in (2.13), these conditions could, in fact, be relaxed to steepness
in half of the directions, i.e., in $d$-dimensions only. 

We use Lemma~2.3 to prove the Main Lemma, where we will set $n=2d$ and $\rho=2^{2d}+1$. 
We first state a direct corollary of Proposition~2.2 in \cite{B3}, by assuming $u_N$ and $E_N$ are {\it fixed}, instead 
of varying with $N$. (The ``$E$" in \cite{B3} is set to be $0$ here.)

\proclaim{Lemma~2.4} Let $I$ be an interval in $(1/2, 3/2)^{2d}$ as in the Main Lemma, and $u_N=u^{(0)}$, $E_N=E^{(1)}$ (in (2.8)) for all $N$. 
There exists $N_0=N_0(a)$, such that for all $N\geq N_0$, there exists $\Cal A_N\subset I$
satisfying (2.20), and $$\text{meas }(I\backslash \cap_{N\geq N_0}\Cal A_N)\to 0,$$ as $a\to 0$.
On $\Cal A_N$, (2.22)-(2.24) hold.
\endproclaim

\demo {Proof}
Choose $N_0=|\log a|^s$ for some $s>1$. Set 
$$D=\sum_{k=2}^d(j_k\cdot\lambda_k+\theta_k)^2-(\tilde j\cdot\lambda)^2.\tag 2.33$$
From (2.7) and (2.10), to prove (2.23) and (2.24) at $N=N_0$, it suffices that
$$|(j_1\cdot\lambda_1+\theta_1)^2+D|\geq a^{p+1}\tag 2.34$$
for all $j=(j_1, j_2, ..., j_d)\in [-N_0, N_0]^{2d}$.  This leads to excise a set in $\theta_1$ of measure satisfying (2.22),
if $0<s\tau<1$ and $s\sigma>1$. No excision in $\lambda$ is needed for this step, so (2.20) is trivially satisfied. 

The iteration to larger scales,  $N>N_0$,  uses Lemmas~2.2 and 2.3. In order to import directly the proof of Proposition~2.2 in \cite{B3},
we shall not use the special structure of the $\Bbb Z^{2d}$ action in (2.13). Therefore we double the dimension and introduce
$$\tilde \theta=(\tilde \theta_1, \tilde \theta_2, ..., \tilde \theta_d)
=(\tilde \theta_{1,1}, \tilde \theta_{1,2}, \tilde \theta_{2,1}, \tilde \theta_{2,2}, ..., \tilde \theta_{d,1}, \tilde \theta_{d,2})\in\Bbb R^{2d}.\tag 2.35$$
The covariance of the $\Bbb Z^{2d}$ action on $\Bbb R^{2d}$ of the corresponding linearized operator $F'(\tilde\theta)$ is then: 
$$(\tilde \theta_{1,1}, \tilde \theta_{1,2}, \tilde \theta_{2,1}, \tilde \theta_{2,2}, ..., \tilde \theta_{d,1}, \tilde \theta_{d,2})
\mapsto (\tilde\theta_{1,1}+\tilde\theta_{1,2}+j_1\cdot \lambda_1, \tilde\theta_{2,1}+\tilde\theta_{2,2}+j_2\cdot \lambda_2, ..., \tilde\theta_{d,1}+\tilde\theta_{d,2}+j_d\cdot \lambda_d).\tag 2.36$$ 
We note that the right hand side is  {\it independent} of 
$$\theta^{-}_i=\tilde \theta_{i,1}-\tilde \theta_{i,2},\text{ for all } i=1, 2, ..., d.$$ 
It follows that 
$$\theta_i=\tilde \theta_{i,1}+\tilde \theta_{i,2},\, i=1, 2, ..., d,\tag 2.37$$
to return to (2.13). 

Let $\tilde \Theta_N\subset\Bbb R^{2d}$ be the set, on the complement of which, (2.23) and (2.24) (in the argument $\tilde\theta$) hold.
From the above discussion, the set 
$\tilde \Theta_N$ is {\it independent} of 
$\theta^{-}_i$, $i=1, 2, ..., d$.
Clearly, for the initial estimate at scale $N_0$, we may proceed as in (2.33)-(2.34) and obtain 
$\tilde \Theta_{N_0}$, satisfying 
$$\text{meas } [\tilde\theta_{i, 1}| \forall \text { fixed }\tilde\theta_{i, 2}, \tilde\theta_k, k\neq i; \tilde\theta\in \tilde\Theta_{N_0}]\leq e^{-N_0^\tau} \, (\tau>0),$$
$$\text{meas } [\tilde\theta_{i, 2}| \forall \text { fixed }\tilde\theta_{i, 1}, \tilde\theta_k, k\neq i; \tilde\theta\in \tilde\Theta_{N_0}]\leq e^{-N_0^\tau} \, (\tau>0),$$
for all $i=1, 2, ..., d.$
Since $a$ is fixed, $N_0$ is fixed; $u^{(0)}$ and $E^{(0)}$ are fixed, $T$ is a fixed operator; using $\tilde\Theta_{N_0}$, we are in exactly the same setting
as \cite{B3} and \cite{JLS} in $2d$ dimensions.

Proposition~2.2 \cite{B3} and Theorem~4.1 \cite{JLS} are directly applicable. With $2d$ replacing $d$, we obtain that for all scales $N\geq N_0$, there is a
good frequency set $\Cal A_N$, $\text{deg } \Cal A_N\leq N^{8d}$, such that the conclusions in (2.22)-(2.24) hold, with $\tilde\Theta_N$ replacing $\Theta_N$, and
$$\text{meas } [\tilde\theta_{i, 1}| \forall \text { fixed }\tilde\theta_{i, 2}, \tilde\theta_k, k\neq i; \tilde\theta\in \tilde\Theta_N]\leq e^{-N^\tau} \, (\tau>0),\tag 2.38$$
$$\text{meas } [\tilde\theta_{i, 2}| \forall \text { fixed }\tilde\theta_{i, 1}, \tilde\theta_k, k\neq i; \tilde\theta\in \tilde\Theta_N]\leq e^{-N^\tau} \, (\tau>0),\tag 2.39$$
for all $i=1, 2, ..., d.$ 
Using (2.37), (2.38)-(2.39) become 
$$\text{meas } [\theta_{i}| \forall \text { fixed } \tilde\theta_k, k\neq i; \tilde\theta\in \tilde\Theta_N]\leq e^{-N^\tau} \, (\tau>0).$$
Since $\tilde\Theta_N$ is independent of $$\theta^{-}_i=\tilde \theta_{i,1}-\tilde \theta_{i,2},\text{ for all } i=1, 2, ..., d,$$
this with (2.37) lead to (2.22). The measure estimate in the Lemma also follows. 
\hfill $\square$
\enddemo

\demo{Proof of the Main Lemma} Scale $N=N_0$ is already proved in Lemma~2.4. 
To obtain (2.23) and (2.24) at larger scales, we use induction. 
Below we keep to the \cite{JLS} notations. 

The induction uses $3$ scales: $N_1$, $N_2=N_1^{2/c_1}$ and $N_3=e^{N_1^{c_1}}$, where $c_1>0$,
and $c_1= \tau$ here. Assume that the Main Lemma holds at scales $N_1$ and $N_2$, we shall show that it holds at $N_3$.

At scale $N_1$, if $\lambda\in\Cal A_{N_1}$, then for $\tilde\theta\notin\tilde\Theta_{N_1}$:
$$\Vert [T_{N_1}(u_{N_1}, E_{N_1})(\tilde\theta)]^{-1}\Vert_{\text{Op}}\leq e^{N_1^\sigma}\, (1>\sigma>\tau>0),\tag 2.40$$
$$|[T_{N_1}({u_{N_1}, E_{N_1}})(\tilde\theta)]^{-1}(j,j')|\leq e^{-\beta|j-j'|} (0<\beta<\gamma),\,\forall |j-j'|>N_1/10.\tag 2.41$$
Due to geometric reasons in the induction (paving) process, aside from cubes: 
$$Q_N=[-N, N]^{2d},\tag 2.42$$ we also need to consider
regions of the form:
$$Q_N=[-N, N]^{2d}\backslash \{n\in\Bbb Z^{2d}:n_i\zeta_i 0, 1\leq i\leq 2d\},\tag 2.43$$
where for $i=1, 2, ..., 2d$, $\zeta_i\in\{<, >, \emptyset\}^{2d}$ and at least two $\zeta_i$ are not $\emptyset$. 
We assume that (2.40) and (2.41) hold for $Q_{N_1}$ as well. 

Theorem~2.7 in \cite{JLS} is applicable and gives the following:
\noindent There is a semi-algebraic set $\Cal A_3\subset\Cal A_{N_1}$, with $\text{deg }\Cal A_3\leq N_3^{8d}$ and 
$$\text{meas }(\Cal A_{N_1}\backslash \Cal A_3)\leq N_3^{-c_3},\tag 2.44$$
with $c_3=8dc_1=8d\tau$, such that if $\lambda\in\Cal A_3$, then for any $\tilde {\vartheta}\in \Bbb R^{2d}$, there 
exists $\tilde N\in [N_3^{c_3}, N_3^{c_4}]$ and annulus 
$$\Gamma:=[-\tilde N, \tilde N]^{2d}\backslash [-\tilde N^{\frac{1}{20d}}, \tilde N^{\frac{1}{20d}}]^{2d},$$
such that for all $k\in \Gamma$,  $Q_{N_1}(u_{N_1}, E_{N_1})(\tilde\vartheta)+k$ satisfy (2.40) and (2.41). Using (2.15), (2.16),
(2.18) and (2.19) between scales $N_1$ and $N\geq N_3$, we obtain 
$$\Vert [T_{N_1}(u_{N}, E_{N})(\tilde\vartheta)]^{-1}\Vert_{\text{Op}}\leq (1+e^{-\gamma N_1}) e^{N_1^\sigma}\, (1>\sigma>\tau>0),\tag 2.45$$
$$|[T_{N_1}({u_{N}, E_{N}})(\tilde\vartheta)]^{-1}(j,j')|\leq (1+e^{-\gamma N_1})e^{-\beta|j-j'|}\, (0<\beta<\gamma),\,\forall |j-j'|>N_1/10,\tag 2.46$$
where $T_{N_1}$ now denotes restrictions to $Q_{N_1}+k$, $k\in \Gamma$. 

We are now at the point to apply Theorem~3.6 in \cite{JLS}, which uses Cartan Theorem. For this purpose, we need the estimates (2.23) and (2.24)
at scale $N_2$. Using (2.15), (2.16), (2.18) and (2.19) between $N_2$ and $N\geq N_3$, yields: 
$$\Vert [T_{N_2}(u_{N}, E_{N})(\tilde\theta)]^{-1}\Vert_{\text{Op}}\leq (1+e^{-\gamma N_2}) e^{N_2^\sigma}\, (1>\sigma>\tau>0),$$
$$|[T_{N_2}({u_{N}, E_{N}})(\tilde\theta)]^{-1}(j,j')|\leq (1+e^{-\gamma N_2}) e^{-\beta|j-j'|}\, (0<\beta<\gamma), \forall |j-j'|>N_2/10.$$

Applying Theorem~3.6 \cite{JLS}, we then obtain that for all $N\in[N_3, N_3^2]$, (2.20), (2.22)-(2.24) hold at scale $N$, 
with $$\beta_N=\beta-\Cal O(1)/N_1^\sigma$$
replacing $\beta$ in (2.24), provided
$$\lambda\in \Cal A_N:=\Cal A_3\cap \Cal A_{N_2}.\tag 2.47$$ 

We take an interval of initial scales $N_0\in [(\log|\log a|)^{1/\tau}, |\log a|^s]$,
$0<s\tau<1$. Clearly (2.22)-(2.24) hold for all such $N_0$, and $\Cal A_{N_0}=I$. 
For example, for the scale $(\log|\log a|)^{1/\tau}$, modifying (2.34) to require
$$|(j_1\cdot\lambda_1+\theta_1)^2+D|\geq 2|\log a|^{-1},$$
for all $j=(j_1, j_2, ..., j_d)\in [-N_0, N_0]^{2d}$, leads to the desired estimates.
Denote the sub-exponential induction in scales by the function $f$: e.g., $N_3=f(N_1)=e^{N_1^\tau}$, $\tau>0$, since $f^2(x)>f(x+1)$,
e.g., $N_3^2=e^{2N_1^\tau}> e^{(N_1+1)^\tau}$, the iterates therefore generate all possible scales. Consequently, we obtain that (2.22)-(2.24) hold for all $N$
with exponential rate of decay $\beta_N$ bounded below by
$$\beta_\infty=\beta_{N_0}-\sum_{i=0}^\infty\frac{\Cal O(1)}{[f^{(i)}(N_0)]^\sigma}>\frac{\beta_{N_0}}{2},\tag 2.48$$
where  $f^{(i)}$ is the $i$th iterate of $f$. Set $$\beta=\beta_\infty.\tag 2.49$$

We are only left to prove (2.21). From construction, if both $N-1$ and $N\in [N_3, N_3^2]$, for some $N_3$, then $\Cal A_{N-1}\backslash\Cal A_N=\emptyset$,
$\text{meas }\Cal A_{N-1}\backslash\Cal A_N=0$. Otherwise, if $N= N_3^2+1$, then from (2.44) one needs
to make an additional excision of measure less than
$$(N_3^2+1)^{-c_3}=N^{-8d\tau}$$ from 
the set $\Cal A_{N_1'}$, satisfying the inclusions,
$$\Cal A_{N_1'}\supset\Cal A_{N_2'}\supset \Cal A_{N_3},\tag 2.50$$ 
where $N_1'=(\log N)^{1/\tau}\ll N_3$ and $N_2'={N_1'}^{2/\tau}\ll N_3$, are the two scales in the induction to reach scale $N$. Since $\Cal A_{N_3}=\Cal A_{N-1}$,
the inclusion in (2.50) gives (2.21) with $c=8d\tau$. \hfill$\square$   
\enddemo

\subheading{2.5. How to use the Main Lemma} 

We anticipate in the next few lines the application of the Main Lemma to the nonlinear analysis in sect.~3.

From the Newton scheme (2.9), and its multiscale realizations, the $u_N$, hence $T_N$, in the Main Lemma are defined only
on intervals $I$ such that an appropriate restricted $F'=T$ is invertible. (This is one of the main differences
with linear theory, where the operators are given, and therefore defined a priori on all of the parameter space.)  
So in the nonlinear application in sect.~3, the intervals $``I"$ in the Main Lemma, will vary with $N$, and 
we shall apply the Lemma to {\it each} interval $I$ in $\lambda$, on which $T_N$ is defined. Note that the
measure estimate in (2.21) is {\it per} interval. To control the total excised measure, the nested properties 
of different generations of intervals, which we already had a glimpse of in the proof of the Main Lemma, shall come into play. 

\head{\bf 3. Nonlinear construction -- proof of the Theorem }\endhead

Our goal is to seek space quasi-periodic solutions in the form (QP) to the NLS (1.1).
Recall that it leads to the nonlinear matrix equation (2.1), which we denote by $F(u)=0$ in (2.3). 
(Recall also that $u$ now stands for $\hat u$.)
The equations are divided into the $Q$-equations in (2.4), leading to (2.7), which yields $E$, 
$$E=(\tilde j\cdot\lambda)^2-(2^d/a)(u)^{*2p+1}|_{\tilde j},\tag 3.0$$
and the $P$-equations (2.5), which are used to solve for $u$.
The initial approximation $u^{(0)}$ is given by $\tilde u$ (1.7). Substituting $u^{(0)}$ into (3.0)
gives 
$$E^{(1)}=(\tilde j\cdot\lambda)^2-(2^d/a)(u^{(0)})^{*2p+1}|_{\tilde j}.$$

The $P$-equations are solved using a Newton scheme and iteration in scales.
Let $M$ be a large integer, and $M^r$, $r=1, 2, ...$, the geometric sequence of scales.
Denote by $u^{(r)}$ the $r$th approximation, and the increment $$\Delta u^{(r)}=u^{(r)}-u^{(r-1)}.$$
We define 
$$\Delta u^{(r)}= -[F'_N(u^{(r-1)}, E^{(r)})]^{-1}F(u^{(r-1)}),\tag 3.1$$
where $N=M^{r}$ and $F'_N(u^{(r-1)}, E^{(r)})$ is the restriction of the linearized operator,
$$F'_N(u^{(r-1)}, E^{(r)})=\text{diag }(\sum_{k=1}^d(j_k\cdot\lambda_k)^2-E^{(r)})-(2p+1)({u^{(r-1)}})^{*2p}*,$$
to the cube $[-N, N]^{2d}$, and where
$$E^{(r)}=(\tilde j\cdot\lambda)^2-(2^d/a)(u^{(r-1)})^{*2p+1}|_{\tilde j}.\tag 3.2$$

Equations (3.2) and (3.1) together with $u^{(0)}$ iteratively solve the $Q$ and the $P$-equations, provided
(3.1) is well-defined for all $r$ and the resulting series converges. The Main Lemma is pivotal in
estimating $[F'_N(u^{(r-1)}, E^{(r)})]^{-1}$ in (3.1), which ensures double exponential convergence of the Newton scheme. We first lay down the induction hypothesis.

Let $M$ be a large positive integer. As earlier, one may assume $a>0$. It consists in showing that the following are satisfied for all $r>0$ and fixed small $a$:

On the {\it entire} $\lambda$ space, namely $(1/2, 3/2)^{2d}$:
\item{(Hi)} $\text{supp }u^{(r)}\subseteq B(0, M^r)$ ($\text{supp }u^{(0)}\subset B(0, M)$).
\item{(Hii)} $\Vert \Delta u^{(r)}\Vert<\delta_r$,
$\Vert \partial \Delta u^{(r)}\Vert<\tilde\delta_r$ with $\delta_{r+1}\ll\delta_{r}$ and $\tilde\delta_{r+1}\ll\tilde\delta_{r}$,
where $\partial$ denotes $\partial_\lambda$ and $\Vert\,\Vert:=\sup_{\lambda}\Vert\,\Vert_{\ell^2(\Bbb Z^{2d})}$.
\item{(Hiii)} $|u^{(r)}(j)|<ae^{-\alpha|j|}\, (\alpha>0)$. 
\bigskip
Using (3.6) and (Hi-iii), the nonlinear eigenvalue
$E^{(r)}$ is $C^1$ in $\lambda$ on $(1/2, 3/2)^{2d}$.
Moreover  by (Hii),
$$| E^{(r)}- E^{(r-1)} |\lesssim \Vert u^{(r)}-u^{(r-1)}\Vert <\delta_r,$$
so that $E^{(r-1)}$ is a $\delta_r$ approximation of $E^{(r)}$. 
\bigskip
Below we continue with the assumptions on the {\it restricted} intervals in $\lambda$ on $(1/2, 3/2)^{2d}$, where one could construct approximate solutions.

\item{(Hiv)} There is a collection $\Lambda_r$ of intervals of size $a^{p+2}M^{-r^C}$, $C>7$, such that 
\item{(a)} On $I\in\Lambda_r$, $u^{(r)}(\lambda)$ is given by a rational function in $\lambda$ of degree at most 
$M^{Cr^3}$. (Consequently, $E^{(r)}$ is rational of degree at most $M^{(2p+1)Cr^3}$ from (3.2).)
\item{(b)} For $\lambda\in\bigcup_{I\in\Lambda_r} I$,

$\Vert F(u^{(r)})\Vert<\kappa_r$,
$\Vert \partial F(u^{(r)})\Vert<\tilde\kappa_r$ 
with $\kappa_{r+1}\ll\kappa_{r}$ and $\tilde\kappa_{r+1}\ll\tilde\kappa_{r}$
\item{(c)} Let $N=M^r$. For $\lambda\in\bigcup_{I\in\Lambda_r} I$, $T=T(u^{(r-1)}):=F'(u^{(r-1)})$ satisfies

$\Vert T_N^{-1}\Vert <a^{-(p+2)}M^{r^C}$,

$|T_N^{-1}(j, j')|<a^{-(p+2)}e^{-\alpha|j-j'|}$, for $|j-j'|>r^{C}$,

where $T_N$ is $T$ restricted to $[-N, N]^{2d}$.
\item{(d)} Each $I\in\Lambda_r$ is contained in an interval $I'\in\Lambda_{r-1}$ and 
$$\text{meas}(\bigcup_{I'\in\Lambda_{r-1}} I'\backslash \bigcup_{I\in\Lambda_{r}}I)<a^{p/5}r^{-5}.$$

The iteration holds with 
$$\delta_r<a^pM^{-(\frac{4}{3})^r}, \, \tilde\delta_r<a^pM^{-\frac{1}{2}(\frac{4}{3})^r}; \kappa_r<a^{2p} M^{-(\frac{4}{3})^{r+2}}, \, \tilde\kappa_r<a^{2p} M^{-\frac{1}{2}(\frac{4}{3})^{r+2}}.\tag W$$

\smallskip 
We remark that the approximate solutions $u^{(r)}$ are defined, a priori, on $\Lambda_r$, but using the derivative estimates 
in (Hii) together with (W), as $C^1$ {\it functions} they can be, and are 
extended to $(1/2, 3/2)^{2d}$, by using a standard argument.
\bigskip

\subheading{3.1. About the induction hypothesis}

Let us provide some intuitions to the hypothesis in (Hiv); (Hi-iii) follow by the construction defined in (3.1), (3.2) and the Newton scheme. 

First of all, since $u^{(0)}$ is independent of $\lambda$ and $E^{(0)}=\tilde E$ in (1.9) is quadratic in $\lambda$, 
$E^{(r)}$ defined in (3.2) and $u^{(r)}(\lambda)$ defined using (3.1), are (formally) clearly rational functions, for $r=1, 2, ...$, which provide the basis for our analysis.  

\item {$-$} Size of the intervals: controlled by the bounds in (Hiv, c), as one may perturb $\lambda$ and retain essentially the same bound. 
           
\item {$-$} Number of intervals: given by the inverse of the size of the intervals. Note that it is obtained by {\it analytic} arguments, and not topological ones.         
            
\item {$-$} Pointwise estimates in (Hiv, c): available at scales $r^C\ll N=M^r$, this means that the Main Lemma is used at {\it much smaller}
scales $r^C$, cf. (2.24).   

\item {$-$} The intervals in $\Lambda_r$: the smaller the $r$, the smaller number of intervals, hence complexity; moreover there is the nested 
property exhibited in (Hiv, d). This will be essential when applying the Main Lemma.

\item {$-$} The induction: consists of the initial steps and the general steps. The initial steps are direct perturbations using small amplitude $a$;
the general steps use the Main Lemma and then convert the $\theta$-estimates into $\lambda$ estimates, leading to successive generations
of $\Lambda_r$.       







\subheading{3.2. The initial steps: $r\leq R$} 

We start with the initial steps.
In the Lemma below, for simplicity, $T_N$, $N=M^r$, stand for $T_N(u^{(r-1)}, E^{(r)})$.

\proclaim{Lemma 3.1} There is a set $\Cal B_N$ in $\lambda$, with
$\text{meas } \Cal B_N<a^{p/5}$, such that on $(1/2, 3/2)^{2d}\backslash\Cal B_N$, 
$$\aligned &\Vert T_N^{-1}\Vert <a^{-(p/2)},\\
&|T_N^{-1}(j, j')|<a^{-(p/2)}e^{-|\log a||j-j'|},\endaligned\tag 3.3$$
for all $N\leq e^{|\log a|^{5/6}}$. 
\endproclaim
\demo{Proof} This follows from perturbation of the diagonals. From (3.2),
$$E^{(1)}=\sum_{k=1}^d(\tilde j_k\cdot\lambda)^2+\Cal O(a^{2p})=(\tilde j\cdot\lambda)^2+\Cal O(a^{2p}),$$
it suffices if 
$$|\sum_{k=1}^d(j_k\cdot\lambda_k)^2-(\tilde j\cdot\lambda)^2|>2a^{p/2},\tag 3.4$$
for all $j=(j_1, ...,  j_k, ..., j_d)\in [-N, N]^{2d}\backslash\Cal S$. For each $j\in [-N, N]^{2d}\backslash\Cal S$, 
it is easy to see that 
$$\aligned &\sum_{k=1}^d[(j_k\cdot\lambda_k)^2]-(\tilde j\cdot\lambda)^2
=\sum_{k=1}^d[(j_k\cdot\lambda_k)^2-(\tilde j_k\cdot\lambda_k)^2]\\
&=\sum_{k=1}^d[(j_k-\tilde j_k)\cdot\lambda_k][(j_k+\tilde j_k)\cdot\lambda_k ] \not\equiv 0,\endaligned$$
by setting $\lambda= (1, 0, ..., 0)$ and the Diophantine condition (D).
It is a quadratic polynomial in $\lambda$. Summing over $j$ then gives the measure estimate for (3.4) to hold. The   
norm estimate in (3.3) follows from (3.4) by simple perturbation; while the pointwise estimate by resolvent series expansion.
\hfill$\square$
\enddemo

\proclaim{Corollary} Set $R=|\log a|^{3/4}$, (Hi-iv) and (W) hold for $1\leq r\leq R$.
\endproclaim

\demo{Proof}
We first address (Hiv). Using Lemma~3.1 for the first $R$, $R=|\log a|^{3/4}$, 
steps of the induction, (Hiv, c) is verified, with 
$$\alpha=\Cal O(|\log a|),\tag 3.5$$
for all scales $N$, $$N\in [M,  M^{|\log a|^{3/4}}],$$ 
with corresponding sets of intervals $\Lambda_r$, $r=(\log N/\log M)\leq R$:
$$\bigcup_{I\in\Lambda_r} I \subseteq (1/2, 3/2)^{2d}\backslash\Cal B_N.$$
The nested property (Hiv, d) is manifest.
On each $I$, (3.3) is satisfied. Clearly (3.4) and hence (3.3) are stable under perturbations of size $a^{p+2}$. So the intervals $I$ are of size $\Cal O(a^{p+2})$.

To prove (Hiv, a), we use induction. When $r=0$, $u^{(0)}$ is independent of $\lambda$. Assume that it holds at stage $r$, 
$$\text{deg } u^{(r)}\leq M^{Cr^3}.$$
Appealing to the definition (3.1) and using the expression below it, we obtain
$$\text{deg } u^{(r+1)}\lesssim 2p\text{ deg } u^{(r)} M^{2d(r+1)}<M^{C(r+1)^3},$$
where the volume factor $M^{2d(r+1)}$ stems from the determinant used to compute the inverse.
The above argument evidently holds for all $r=1, 2, ...$.

We are left with (Hiv, b). When $r=0$, 
$$F(u^{(0)})=\Cal O(a^{2p}),$$ 
and 
$$\partial F(u^{(0)})=0.$$
On $\Lambda_r$, $\Delta u^{(r)}$ is constructed using (3.1): 
$$\Delta u^{(r)}= -[F'_N(u^{(r-1)})]^{-1}F(u^{(r-1)}).$$
Using (3.3), this gives, when $r=1$,
$$\Vert \Delta u^{(1)}\Vert=\Cal O(a^{3p/2}),$$
more over (Hiii) is satisfied with $\alpha$ satisfying (3.5).
Similarly 
$$\Vert \partial \Delta u^{(1)}\Vert=\Cal O(a^{p}).$$
So the first two expressions in (W) at $r=1$ are satisfied. 

To verify the other two expressions, we write
$$\aligned F(u^{(1)})&= F(u^{(0)})+F'(u^{(0)}) \Delta u^{(1)}+\Cal O ((\Delta u^{(1)})^2)\\
&=(T-T_N)\Delta u^{(1)}+\Cal O ((\Delta u^{(1)})^2)\\
&=-[(T-T_N)T^{-1}_N]F(u^{(0)})+\Cal O (\Vert T^{-1}_N\Vert^2F(u^{(0)})^2)\\
&<a^{2p+1},\endaligned$$
using (Hi-iii) at $r=1$, (3.1) and (3.3). This verifies the third expression in (W) at $r=1$. Similarly we may verify
the fourth expression.

The extension argument in sect.~10, (10.33-10.37) in \cite{B1}, then proves (Hi-iii) on the entire $\lambda$ space, moreover (W) is satisfied at $r=1$.
Iterating the above arguments, we prove (Hi-iv) and (W) for all $r\leq R$ with $\alpha$ satisfying (3.5).  
(For details of the iteration to prove (W), see Lemma~5.2 and its proof in \cite{W1}. This is rather
routine and clearly holds for all $r=1, 2, ...$)\hfill$\square$ 
\enddemo

\noindent{\it Remark.} As a side, we mention that the semi-algebraic sets $\Cal B_N$ can be described by the violation of $(2N+1)^{2d}$ quadratic polynomial 
inequalities in (3.4), Lemma~2.1 gives that the number of connected components in $\Cal B_N$ is
bounded above by $\Cal O(N^{4d^2})$. The set $\Cal B_N$ is, moreover, independent of $u^{(0)}$, in fact, all $u^{(k)}$ for $k\leq R-1$. 
\bigskip
Unlike the first $R$ steps, however, the iterations to subsequent scales use the $\theta$ estimates
and the Main Lemma to make excisions in $\lambda$, in order to fulfill (Hiv,~c,~d). Afterwards, the same induction arguments
used in the first $R$ steps, will validate
(Hi-iii, Hiv,~a,~b) and (W) for all $r>R$.

\subheading{3.3. The general steps: $r\geq R$}

Assume (Hi-iv) hold at stage $r$. To construct $u^{(r+1)}$, we need to control 
$$T_N^{-1}(u^{(r)}) \text { with }N=M^{r+1}.$$
This requires another excision in $\lambda$, which will lead to the next set of intervals $\Lambda_{r+1}$. 

To simplify notations, given two sets of intervals $Z_1$ and $Z_2$, we say that 
$$Z_2\subset Z_1,$$
if for all $I\in Z_2$, there exists $I'\in Z_1$, such that
$I\subset I'$.  We also define
$$\text{meas } (Z_1\backslash Z_2)=\text {meas }(\bigcup_{I'\in Z_1} I'\backslash\bigcup_{I\in Z_2} I.)$$

Cover $[-M^{r+1}, M^{r+1}]^{2d}$ by $[-M^{r}, M^{r}]^{2d}$ and smaller
cubes $[-M_0, M_0]^{2d}+J$, with $M^r/2<|J|<M^{r+1}$ and $M_0\ll N$ to
be specified shortly. Let $r\geq R$. For simplicity, we drop the prefactors in $a$ in (Hi-iv), since
they are fixed, and only keep track of variations in $r$. 
\smallskip 
\noindent $\bullet$ {\bf $\theta$-estimates of} $T^{-1}_{M_0}(\theta)$

The following Lemma provides $\theta$ estimates on the $M_0$-cube centered at the origin. 
\proclaim{Lemma 3.2} Assume (Hi-iv) hold at stage $r$, scale $N=M^r$. Set 
$$M_0=(\log N)^C=r^C(\log M)^C,\tag 3.6$$ with $C>7/c$, and $c$ as in (2.21);
define $$r_0=\frac{\log M_0}{\log M},\tag 3.7$$ and 
$$\tilde r_0:=r_0\frac{\log M}{\log 4/3}<2C\frac{\log r}{\log 4/3}\ll r.\tag 3.8$$
Then there is  $\Lambda'_{r+1}\subset \Lambda_r$, so that on $\Lambda'_{r+1}$, the following estimates hold:
$$\aligned \Vert T_{M_0}^{-1}(u^{(\tilde r_0)}, E^{(\tilde r_0+1)})(\theta)\Vert &<e^{M_0^\sigma}, \, 0<\sigma<1,\\
|T_{M_0}^{-1}(u^{(\tilde r_0)}, E^{(\tilde r_0+1)})(\theta)(x, y)| &< e^{-\alpha |x-y|},\, \alpha>0,\endaligned\tag 3.9$$
for all $x, y$ such that $|x-y|>M_0/10$, provided $\theta$ is in the complement of a set $\Theta_{M_0}$, whose sectional measures
satisfy
$$\text{meas } [\theta_i| \forall {\text fixed }\theta_k, k\neq i; \theta\in \Theta_N]\leq e^{-M_0^\tau} \, (\tau>0).\tag 3.10$$
\endproclaim

\noindent {\it Remark.} The expressions in (3.6)-(3.8), if not manifestly integers, are understood to be the integer part.

\demo{Proof} 
For the first $r\leq R$ steps, direct perturbation in $a$ proves that (3.9) and (3.10) hold on $\Lambda_r$ with $\alpha=\Cal O(|\log a|)$, 
without additional excisions. (Here one may assume $M_0\geq M$.)

For $r\geq R$, set $$N_1=(\log M_0)^{1/\tau}<(2C\log r)^{1/\tau}$$ from (3.6), and 
$$\tilde r=\frac{\log N_1}{\log M}$$
and 
$$\tilde{\tilde r}=\tilde r \frac{\log M}{\log 4/3}<(10/\tau)\log\log r\ll \tilde r_0\ll r.\tag 3.11$$
(In the language of the Proof of the Main Lemma, $M_0:=N_3$.)
To apply the Main Lemma, fix $I\in \Lambda_{\tilde{\tilde r}}$. 
By the choice of $\tilde{\tilde r}$ and using (W), on $\Lambda_{\tilde{\tilde r}}\cap \Lambda_{\tilde r_0}$, 
$$\Vert T_{N_1}(u^{(\tilde{\tilde r})}, E^{(\tilde{\tilde r}+1)})(\theta)-T_{N_1}(u^{(\tilde r_0)}, E^{(\tilde r_0+1)})(\theta)\Vert \lesssim\delta_{\tilde{\tilde r}}\leq e^{-\tilde\alpha N_1},\, \tilde\alpha>\alpha>0.$$
(Here ``$\cap$" is in the sense of intersections of the intervals in the two sets; note from (Hiv,~d) that each interval in $\Lambda_{\tilde r_0}$ is {\it contained} in an interval 
in $\Lambda_{\tilde{\tilde r}}$.)
There are at most $$M^{(\tilde{\tilde r})^C}\simeq M^{(\log\log r)^C}$$
such intervals in $\Lambda_{\tilde{\tilde r}}$, by using (Hiv) at stage 
$$\tilde{\tilde r}\sim\log\log r\ll r.$$
Denote the intersection over $I$ of the good sets by $\tilde{\Cal A}_{M_0}$ and let $\Lambda'_{r+1}=\Lambda_{r}\cap\tilde{\Cal A}_{M_0}$, then 
$$\text{meas } \Lambda_r\backslash\Lambda'_{r+1}<M^{(\log\log r)^C}/M_0^c<M^{(\log\log r)^C}/r^{Cc}<1/r^6,$$
if $Cc>7$, using (2.21). Here we appealed again to (Hiv,~d), but at stage $r$, namely that each interval in $\Lambda_r$ is contained in an interval in $\Lambda_{\tilde r_0}$. On $\Lambda'_{r+1}$, (3.9)-(3.10) hold.
\hfill$\square$
\enddemo
\smallskip 
\noindent $\bullet$  {\bf Invertibility of} $T^{-1}_{M_0}(\theta=0)$

The projection lemma below, stated as (1.5) in \cite{B3}, converts the $\theta$ estimates in (3.9)-(3.10) for the $M_0$-cube centered at the origin to 
$M_0$-cubes centered at large $J\in\Bbb Z^{2d}$ at $\theta=0$.

\proclaim{Lemma~3.3}
Let $S\subset [0, 1]^{n_1}\times [0, 1]^{n_2}:=[0, 1]^{n}$, be a semi-algebraic set of degree $B$ and $\text{\rm meas}_{n}S <\eta, \log B\ll
\log 1/\eta$.
Denote by $(x, y)\in [0, 1]^{n_1}\times [0, 1]^{n_2}$ the product variable.
Fix $\epsilon>\eta^{1/n}$.
Then there is a decomposition
$$
S =S_1 \bigcup S_2,
$$
with $S_1$ satisfying
$$
\text{\rm meas}_{n_1}(\text{Proj}_x S_1)<B^K\epsilon \quad (K>0),
$$
and $S_2$ the transversality property
$$
\text{\rm meas}_{n_2}(S_2\cap L)< B^K\epsilon^{-1} \eta^{1/n}\quad (K>0),
$$
for any $n_2$-dimensional hyperplane $L$ in  $[0, 1]^{n_1+n_2}$ such that
$$
\max_{1\leq j\leq n_1}|\text{Proj}_L (e_j)|< \frac 1{100}\epsilon,
$$
where $e_j$ are the basis vectors for the $x$-coordinates.
\endproclaim

\noindent{\it Remark.} Lemma~2.3 is, in fact, also derived from Lemma~3.3, cf., the Proof of Lemma~1.20 in \cite{B3}.

\proclaim {Lemma~3.4} 
There exists $\Lambda_{r+1}\subset\Lambda'_{r+1}\subset\Lambda_r$, 
satisfying $$\text{meas } \Lambda_r\backslash \Lambda_{r+1}<1/r^5,$$
provided $C>max (1/\tau, 7/c)$. 
On the intervals in the set $\Lambda_{r+1}$, $T_{[-M_0, M_0]^{2d}+J}^{-1}(u^{(\tilde r_0)})$
satisfy the upper bounds in (3.9) for all $J$ with $M^r/2<|J|<M^{r+1}$.
\endproclaim

\demo{Proof}
We first make estimates on 
$$T_{[-M_0, M_0]^{2d}+J}(u^{(\tilde r_0)}),$$
with $M^r/2<|J|<M^{r+1}$, on the set $\Lambda'_{r+1}$.
Fix $I\in\Lambda_{\tilde r_0}$. 
To apply Lemma~3.3, 
identify the set $\Theta_{M_0}\subset \Bbb R^d$
with the set  
$$\bar\Theta_{M_0}=\Theta_{M_0}\times \{0\}\subset \Bbb R^{2d}.$$
One may assume 
$$\Theta_{M_0}\subset [-4M_0, 4M_0]^d,$$
as otherwise 
$T_{M_0}(\theta)$ is invertible. Make the partition: 
$$\bar\Theta_{M_0}=\cup_K \{[-1, 1]^d+K\}\times \{0\}:=\cup_K \Cal I_K,$$
where $K\in\Bbb Z^d$, satisfying
$$0\leq |K|\leq 4dM_0.\tag 3.12$$

Fix a $K$ and 
let $$\aligned &S_K(\lambda, \bar\theta)=I\times\{\bar\Theta_{M_0}\cap \Cal I_K\}\subset \Bbb R^{4d};\\
&\text{meas }S_K\leq e^{-M_0^\tau}.\endaligned\tag 3.13$$
Below for notational simplicity, we generally write $S$ for $S_K$. 
The set $S$ is described by the opposite of (3.9). 
Replacing the $\ell^2$ norm by the Hilbert-Schmidt norm and since the matrix elements of the inverse 
is the division of two determinants, (3.9) can be expressed as algebraic inequalities
in the matrix elements of degree at most $M_0^C$. Since each matrix element is quadratic in $\bar\theta$ and at most of degree $e^{(C{\log M_0})^3}$ in $\lambda$,
$S$ is of degree at most 
$$\text{deg }S\leq M_0^Ce^{(C{\log M_0})^3}\lesssim e^{({\log M_0})^4}.\tag 3.14$$ 

Let $$J\lambda=(J_1\lambda_1, J_2\lambda_2, ..., J_d\lambda_d)\in\Bbb R^{d},$$ 
where $$J_i\lambda_i=\text{maxv }\{J_{i,1}\lambda_{i,1}, J_{i,2}\lambda_{i,2}\},$$
and ``maxv" denotes maximum in absolute value.
Define the set $\Cal L_I$ to be
$$\Cal L_I=\bigcup_K\bigcup_J S_K(\lambda, \theta=J\lambda),$$
where the union is over $K\in\Bbb Z^d$ satisfying (3.12), and $J\in\Bbb Z^{2d}$, with $M^r/2<|J|\leq M^{r+1}$. 
In the complement of $\Cal L_I$, 
$$T_{[-M_0, M_0]^{2d}+J}$$ satisfies (3.9),
for all $J$ satisfying $M^r/2<|J|\leq M^{r+1}$.

The measure of $\Cal L_I$ is estimated using Lemma~3.3. 
From (3.13) and (3.14), the set $S$ satisfies
$$\aligned \text{deg }S&\lesssim e^{({\log M_0})^4}\lesssim e^{({\log \log N})^4},\\
\text{meas }S&\leq e^{-M_0^\tau}=e^{-(\log N)^{C\tau}}.\endaligned\tag 3.15$$
Since $J$ satisfies
$$M^r/2<|J|\leq M^{r+1},$$
equivalently 
$$N/2<|J|\leq MN,$$ 
we have 
$$\log \text{deg } S\ll \log |J|\ll -\log \text{meas } S.\tag 3.16$$

So the proof in (3.6)-(3.26) in \cite{B3} remains valid,
leading to the conclusion (3.9). Alternatively, one may
view this as a generalization to arbitrary dimensions $d$, of Proposition~5.1 in \cite{BGS}, and its proof.
We therefore have 
$$\text{meas }\Cal L_I\leq N^{-c'}=M^{-c'r},$$
for some $c'=c'(d)>0$, provided $C\tau>1$.
 
Using (Hiv) at stage $\tilde r_0$, the number of intervals at stage $\tilde r_0$ is bounded above by $$M^{{\tilde r_0}^C}\simeq M^{(\log M_0)^C}\simeq M^{(\log r)^C}.$$
Let $$\Cal L=\bigcup_I\Cal L_I,$$ then 
$$\text{meas }\Cal L\leq M^{(\log r)^C}\cdot M^{-c'r}\leq M^{-c'r/2},\, c'>0.\tag 3.17$$

Define $$\Lambda_{r+1}=\Lambda'_{r+1}\backslash\Cal L.$$ 
The set $\Lambda_{r+1}$ satisfis (Hiv, d) at stage $r+1$, adding (3.12) and (3.17) verifies the measure estimate. On the set, $T_{[-M_0, M_0]^{2d}+J}^{-1}(u^{(\tilde r_0)})$
satisfy the upper bounds in (3.9) for all $J$ with $M^r/2<|J|<M^{r+1}$.
\hfill$\square$
\enddemo

\smallskip
\noindent $\bullet$ {\bf Invertibility of} $T^{-1}_{M^{r+1}}(u^{(r)})$

To verify (Hiv, c) at stage $r+1$, we use resolvent expansion and cover $[-N, N]^{2d}=[-M^{r+1}, M^{r+1}]^{2d}$ by the big cube $[-M^r, M^r]^{2d}$ and smaller $M_0$-cubes, and
use the estimates on $T_{M^r}^{-1}(u^{(r)})$ and $T_{[-M_0, M_0]^{2d}+J}^{-1}(u^{(\tilde r_0)})$, $M^r/2<|J|<M^{r+1}$. We have, using (W), 
$$\Vert T_{M^r}(u^{(r)})-T_{M^r}(u^{(r-1)})\Vert\leq\delta_{r-1}\leq M^{-(\frac{4}{3})^{r-1}},$$
and 
$$\Vert T_{[-M_0, M_0]^{2d}+J}(u^{(r)})-T_{[-M_0, M_0]^{2d}+J}(u^{(\tilde r_0)})\Vert\leq \delta_{\tilde r_0}\leq e^{-\tilde\alpha M_0}< M^{-r^C},\tag 3.18$$
for sufficiently large $M$.
(Recall that $\tilde r_0$ is constructed so that the first inequality in (3.18) holds.) 
From (Hiv,~c) at stage~$r$,
$$\Vert T^{-1}_{M^r}(u^{(r-1)})\Vert\leq M^{r^C},\tag 3.19$$
so $$\Vert T^{-1}_{M^r}(u^{(r)})\Vert \leq M^{r^C}+M^{-(\frac{4}{3})^{r-1}}<(1+M^{-(\frac{4}{3})^{r-1}})M^{r^C}.\tag 3.20$$
From Lemma~3.4 and (3.9)
$$\Vert T_{[-M_0, M_0]^{2d}+J}^{-1}(u^{(\tilde r_0)})\Vert\leq e^{M_0^\sigma}\leq M^{r^{\sigma C}}\ll M^{r^C},\tag 3.21$$
so using (3.18), 
$$\Vert T_{[-M_0, M_0]^{2d}+J}^{-1}(u^{(r)})\Vert\leq M^{r^{\sigma C}}+ M^{-r^C}<2M^{r^{\sigma C}}.\tag 3.22$$
The bound in (3.18) shows, moreover, that the pointwise estimates in (3.9) holds
for $T_{[-M_0, M_0]^{2d}+J}^{-1}(u^{(r)})$ with the prefactor $(1+e^{-(\tilde\alpha-\alpha)M_0})$ in lieu of $1$. 

The pointwise estimates on $T^{-1}_{M^r}(u^{(r-1)})$ also essentially hold for $T^{-1}_{M^r}(u^{(r)})$.
This is verified as follows. (Hiv, c) at stage $r$ gives
$$\aligned
\Vert [T_{M^r}(u^{(r-1)})]^{-1}\Vert&\leq M^{r^C},\\
|[T_{M^r}(u^{(r-1)})]^{-1}(k,k')|&\leq e^{-\alpha|k-k'|}\, (|k-k'|>r^C).\endaligned\tag 3.23$$
We write
$$\aligned T_{M^r}(u^{(r)})&=T_{M^r}(u^{(r-1)})+[T_{M^r}(u^{(r)})-T_{M^r}(u^{(r-1)})]\\
&:=A+B\endaligned\tag 3.24$$

To obtain pointwise estimate on $T_{M^r}^{-1}(u^{(r)})$, we use (3.24) and resolvent series. $A^{-1}$
has off-diagonal decay from (3.23), $B$ has off-diagonal decay from (Hiii) at stage $r$. Iterating the
resolvent series and using (3.20), we obtain 
$$|[T_{M^r}(u^{(r)})]^{-1}(k,k')|\leq e^{-\alpha'|k-k'|}\, (|k-k'|>r^C),\tag 3.25$$
with $\alpha'=\alpha-M^{-r\delta'}$ ($\delta'>0$).

Consequently, this leads to
\proclaim{Lemma~3.5} On the set of intervals in $\Lambda_{r+1}$, there are the following estimates: 
$$
\align
&\Vert [T_{M^{r+1}} (u^{(r)})]^{-1}\Vert < M^{2(r+1)d}M^{r^C}\ll M^{(r+1)^C}\\
&|[T_{M^{r+1}}(u^{(r)})]^{-1}(k, k')|< e^{-\bar\alpha |k-k'|} \text{ for } |k-k'|>(r+1)^C,
\endalign
$$
with $\bar\alpha=\alpha-M^{-(r+1)\bar\delta}$, $\bar\delta>0$. 
\endproclaim
\demo{Proof} This is by an application of Lemma~5.1 in \cite{BW}, using (3.20), (3.22), the corresponding pointwise estimates as in (3.9) and (3.25). \hfill$\square$
\enddemo
Lemma~3.5 reproduces the estimates in (Hiv,~c) at stage $r+1$ with a possible lowering of $\alpha$.
Since $u^{(r+1)}$ is constructed using $T^{-1}_{M^{r+1}} (u^{(r)})$, this also represents a lowering
of $\alpha$ in (Hiii). However, similar to (2.48) and (2.49), the ``final" 
$$\alpha:=\alpha_\infty=\Cal O(|\log a|)>0.$$ 
Since the estimates (3.20) and (3.22) are stable under perturbation of size $M^{-(r+1)^C}$, this produces
the next set $\Lambda_{r+1}$ of intervals of size $M^{-(r+1)^C}$. (Hiv, d) follows by construction.
\hfill$\square$

To summarize, we have proved

\proclaim{Lemma~3.6} Assume (Hi-iv) hold at stage $r$, 
There exists $\Lambda_{r+1}\subset\Lambda'_{r+1}\subset\Lambda_r$, 
satisfying $$\text{meas } \Lambda_r\backslash \Lambda_{r+1}<1/r^5,$$
provided $C>max (1/\tau, 7/c)$, such that (Hiv, c, d) hold at stage $r+1$.
\endproclaim

\smallskip 
\noindent$\bullet$ {\bf Construction of} $u^{(r+1)}$

By definition (3.1), 
$$\Delta u^{(r+1)}=-T_{M^{r+1}}^{-1}(u^{(r)})F(u^{(r)}),$$
and $u^{(r+1)}=u^{(r)}+\Delta u^{(r+1)}$. ($E^{(r+1)}$ is already known from (3.2), using $u^{(r)}$.)  Using the same argument as in the proof of the Corollary, verifies (Hiv, a) at stage $r+1$; (Hiv,~b) follows
by direct computation. The derivative estimate in (Hiv, b) leads to the derivative estimate in (Hii).
The extension argument mentioned earlier then shows that
(Hi-iii) hold at stage $r+1$ as well. The induction from
step $r$ to step $r+1$ is thus complete. (Cf. \cite{BW}, sect.~6, (6.1)-(6.20).) \hfill$\square$

\demo{Proof of the Theorem}
The induction process above solves iteratively the $Q$ and the $P$-equations, with the convergence estimates in (W). The measure estimate (1.10)
follows from (Hiv,~d):
$$\text{meas } \Lambda \geq 1-a^{p/5}\sum_{r\geq 1} {r^{-5}}>1-a^{p/6}.$$
The estimates in (1.11) follow from (Hii), (3.2) and (W),
and prove the Theorem. 
\hfill$\square$
\enddemo


\enddocument
\end